\documentclass{article}
 
\usepackage{amssymb,amsmath,amsthm} 
\usepackage{graphicx}
\usepackage[latin1]{inputenc}
\usepackage[a4paper, centering, text={6.5in, 10in}]{geometry} 
\usepackage{authblk}

\newtheorem{theorem}{Theorem}[section]

\newtheorem{proposition}[theorem]{Proposition}
\newtheorem{corollary}[theorem]{Corollary}
\newtheorem{definition}[theorem]{Definition}
\newtheorem{example}{Example}[section]
\newtheorem{remark}[example]{Remark}

\providecommand{\keywords}[1]{\textbf{\textit{Index terms---}} #1}

\begin{titlepage}

\author{Francisco Mota}
\affil{Departamento de Engenharia de Computação e Automação\\
Universidade Federal do Rio Grande do Norte -- Brasil\\
e-mail:mota@dca.ufrn.br}

\date{\today}

\title{Fourier Series and Transforms via Convolution}

\end{titlepage}

\begin{document}

\maketitle

\begin{abstract}
In this paper we show an alternative way of defining Fourier Series and Transform by using the concept of convolution with exponential signals. 
This approach has the advantage of simplifying proofs of transforms properties and, in our view, may be interesting for educational purposes.


\end{abstract}

\section{Introduction}

Fourier Series and Transform \cite{lathi} are pivotal topics in any course of Signals and Systems for engineering. Their 
use is widespread in most engineering courses generally because it help us to solve and/or understand certain 
operations involving signals (e.g. derivation,
integration, translations, etc) that appears in the so-called time-domain  as other operation 
(generally simpler) in another domain denominated 
frequency domain, and vice-versa. Our aim in this note is to present a new formulation for 
Fourier series and transform by exploring
its close connection with another fundamental operation in the context of signal and systems theory that 
is the convolution \cite{lathi} (see also Section~\ref{adconv}). 
The main result of the paper is Proposition~\ref{fsconv} in Section~{\ref{sfourier}, 
which presents another formulation for the Exponential Fourier series. In sections \ref{dtfourier} and \ref{tfourier} we extend
the idea to give a new formulation for the Fourier Transform and Discrete Fourier Transform (DFT), respectively.

\section{Signals and convolution}\label{adconv}

A signal is generally represented as a complex-valued function and which is said to be
{\em analog}  when the domain is the set of real numbers, or {\em discrete} 
when the domain is the set of integers\footnote{The independent variable (domain) may have dimension of time (e.g. seconds) 
or also frequency (e.g. radians/second).}, 
that is:
\begin{eqnarray*}
f:&& \mathbb{R} \to \mathbb{C} \qquad\text{(Analog signal)}\\
  && t \mapsto f(t)\\[0.5cm]
g: && \mathbb{Z} \to \mathbb{C} \qquad\text{(Discrete-time signal)}\\
  && k \mapsto g(k)  
\end{eqnarray*}

As examples we have $f(t)=\cos(\frac{\pi}{2}t)$ as an analog signal and $g(k)=\cos({\frac{\pi}{2}\scriptstyle 10^{-3}}k)$ a discrete signal. 
We can obtain a discrete signal ($f^*$) from an analog signal ($f$) by the process of (periodic) ``sampling", which is mathematically 
implemented as:
\[
f^*(k) = f(kT_s)
\]
where $T_s>0 \in \mathbb{R}$ is denominated ``sampling" interval.\footnote{In practice, the process of sampling is a little more involved, and we
can ``sample" a physical analog signal by using a computer hardware denominated ``Analog-to-Digital converter (or ADC)" \cite{wiki}; 
each sample obtained in this process is a sequence of bits, and so the sampled signal will not only be discrete but also  
digital.} In this situation, we say that the samples of $f$ are spaced in time by an interval $T_s$, and it 
is understood that as $T_s$ tends to zero the discrete signal $f^*$ tends to analog signal $f$, that is $kT_s\to t$ and $f^*(k)\to f(t)$.

Convolution is a binary operation between signals, and we have an analog convolution when both signals involved are analog or a 
discrete convolution when they are discrete signals. We start by defining discrete convolution:

\begin{definition}\label{dconvdef} \em 
The (discrete) convolution between (discrete) signals $f$ and $g$ results in a signal (represented by $f*g$) which is defined as
\begin{equation}\label{dconveq}
(f*g)(k) = \sum_{n=-\infty}^{\infty}f(n)g(k-n).
\end{equation}

\end{definition}

\begin{remark}\label{dconvprop}\em
The infinite (complex) series in Equation~(\ref{dconveq}) is required to be {\em absolutely convergent}, 
in order convolution could share some important properties of other general binary operations, which we present below:

\begin{description}

\item[\em Commutativity:] $f*g = g*f$, for any signals $f$ and $g$. \\
Obs.: Requires infinite series in Equation~(\ref{dconveq}) to be absolutely convergent.

\item[\em Associativity:] $(f*g)*h = f*(g*h)$, for any signals $f$, $g$ and $h$. \\
Obs.: Requires infinite series in Equation~(\ref{dconveq}) to be absolutely convergent.

\item[\em Identity existence:] There exists a signal ``$\delta$", such that $\delta*f = f*\delta=f$, for any signal $f$. 
Signal $\delta$ is defined as
\begin{equation}\label{ddeltadef}
\delta(k) = \begin{cases} 1 &\text{if } k=0\\ 0, &\text{if } k\neq 0\end{cases}
\end{equation}

\end{description}

\end{remark}

We now proceed to define convolution of analog signals (or analog convolution), and as a matter of convenience, we will define it as a 
limit case of discrete convolution. Before all, we introduce the concept of {\em approximated} analog convolution as shown below:

\begin{definition}\label{apconv} \em 
Let be two analog signals $f$ and $g$ and consider their discretization $f^*$ and $g^*$, that is $f^*(k)=f(kT_s)$ and 
$g^*(k)=g(kT_s)$, where $T_s$ is the sampling interval. The approximated (analog) convolution between (analog) signals $f$ and $g$,
results in a signal (represented by $f\tilde{*}g$) which is defined as
\begin{equation}\label{apconveq}
(f\tilde{*}g)(t) = T_s.(f^**g^*)(k) = \sum_{n=-\infty}^{\infty}T_s.f^*(n).g^*(k-n), \quad kT_s\le t < (k+1)T_s
\end{equation}

\end{definition}

\begin{remark} \em 
It is easy to verify that the approximate analog convolution satisfies the same properties for discrete convolution listed in 
Remark~\ref{dconvprop}, but multiplication of discrete convolution formula by the factor $T_s$ requires the identity signal to be
slightly modified; that is, we need to find an analog signal ($\tilde{\delta}$) whose discretization results in discrete signal 
$(1/T_s)\delta$, which is the identity for discrete convolution $T_s(f^**g^*)$. 
While there could be different possibilities, we see that 
\begin{equation}\label{adeltatdef}
\tilde{\delta}(t) = \begin{cases} 1/T_s, &\text{if } -T_s/2 \le t \le T_s/2\\ 0, & \text{otherwise}\end{cases}
\end{equation}
is an analog signal such that its discretization $\tilde{\delta}^*$ results in $(1/T_s)\delta$, as we can see:
\begin{eqnarray*}
\tilde{\delta}^*(k)=\tilde{\delta}(kT_s) & = & \begin{cases} 1/T_s, &\text{if } k=0\\ 0, & \text{if }k\neq 0\end{cases}\\
                    & = & \frac{1}{T_s}\delta(k).
\end{eqnarray*}
And so, we have that $\tilde{\delta}$ defined in Equation~(\ref{adeltatdef}) is an identity signal for the approximated analog convolution.
\end{remark}

We define the (exact) analog convolution just by taking $T_s \to 0$ in Equation~(\ref{apconveq}), and its easy to notice
in this situation that when 
$T_s$ is an infinitesimal ($d\tau$) we have $kT_s\to t$, $nT_s\to \tau$, $f^*(n)\to f(\tau)$, $g^*(k-n)\to g(t-\tau)$ 
and the summand in Equation~(\ref{apconveq}) converges to an (Riemann) integral. 
So that we have:

\begin{definition}\label{aconv}\em 
The (analog) convolution of two (analog) signals $f$ and $g$ is the limit when $T_s\to 0$ of the approximated  convolution (see  
Definition~\ref{apconv}), and it results in a signal $f*g$ defined as:
\begin{equation}\label{aconveq}
(f*g)(t) = \int_{-\infty}^{\infty}f(\tau)g(t-\tau)d\tau
\end{equation} 
\end{definition}

\begin{remark}\em 
In order the analog convolution to be well defined we require that integral in Equation~(\ref{aconv}) to be absolutely convergent, and
under this condition we also can easily prove that, similarly to discrete convolution, analog convolution is  a commutative and associative
binary operation; but we have an issue related to the existence of the identity signal, since when $T_s\to 0$ in Equation~(\ref{adeltatdef}) 
we have that signal $\tilde{\delta}$ becomes undefined at $t=0$. In fact, it is well known that the identity for the analog convolution is 
not a signal (defined as a function), and it is in fact a distribution \cite{distrib}. We just accept it exists as a ``special signal" which is 
the limit of signal $\tilde{\delta}$ (defined in Equation~(\ref{adeltatdef})) when $T_s\to 0$. It is also represented by ``$\delta$", and
so $\delta*f = f*\delta = f$ for any analog signal $f$.

\end{remark}

\subsection{Periodic Signals and Periodic Convolution}\label{pconv}

A periodic (analog) signal $f$ has the property that exists a real number $T>0$ such that $f(t+T) = f(t)$ for all $t\in\mathbb{R}$, 
and similarly, for a discrete signal $g$ to be periodic, it must have be an integer $N>0$ such that 
$g(k+N) = g(k)$ for all $k\in\mathbb{Z}$. With periodic signals,\footnote{We may consider a constant signal as being periodic, where the period 
is any positive value. In analog case, constant signals has no minimum value for the period $T$, while in discrete case the minimum 
value for the  period is $N=1$.} 
it is common to modify the definition of convolution, 
as presented before, in order the interval of integration (or summation) to be reduced to one period of the signal 
(as opposed to the whole domain),\footnote{In fact, the (regular) convolution between periodic signals may diverge due to the fact that periodic
signals are not absolutely integrable (analog) or absolutely summable (discrete).} and then we have the concept of {\em periodic convolution:}

\begin{definition}\label{pconvdef}\em 
The periodic convolution between signals $f$ and $g$, both with same period, results in a periodic signal (with same period of $f$ and $g$), 
represented by $f\circledast g$, and which is defined by:

\begin{eqnarray}
(f\circledast g)(t)  & = & \int_{-T/2}^{T/2}f(\tau)g(t-\tau)d\tau, \qquad f \text{ and } g \text{ are analog signals with same period } T\\
(f\circledast g)(k)  & = & \sum_{n=0}^{N-1}f(n)g(k-n), \qquad f \text{ and } g \text{ are discrete signals with same period } N
\end{eqnarray}

\end{definition}

Periodic convolution can be turned into  a (regular) convolution when one of the periodic signals is switched by its aperiodic component, 
that is, another signal that corresponds just to one period of it and null otherwise:
\[
f\circledast g = f_c*g = f*g_c,
\]
where $f_c$ and $g_c$ are non-periodic signals that corresponds to one period of $f$ and $g$ respectively, and are null otherwise.

\begin{remark}\em 
The convolution between a non-periodic signal $h$ and a periodic signal $f$ results in a signal ($h*f$) which is periodic with 
same period of $f$, so we can mix convolution with periodic convolution, and we have the following associative property 
(in analog or discrete context):
\begin{equation} \label{circass}
(h*f)\circledast g = h*(f\circledast g)
\end{equation}
where $h$ is a non-periodic signal and $f$ and $g$ are both periodic signals with same period.

\end{remark}

\subsection{Some results and properties of convolution}\label{basics}

The most important result, for our purposes, regarding convolution is a very simple fact about convolution with exponential signals:\\[0.2cm]
\fbox{%
\parbox{\textwidth}{%
\em The convolution of an exponential signal with any other signal results in the same exponential signal multiplied by a 
constant factor.
}%
}
\mbox{}\\
We make this statement more precise below:

\begin{proposition}\label{convexp}\em \mbox{}\\

\begin{description}

\item[(a) Analog Convolution with exponential:] Let be $f$ an analog signal and consider $g(t) = e^{at}$, with $a\neq 0\in\mathbb{C}$. Then 
\begin{equation}\label{convexpa}
(f*g)(t) = F(a)g(t)
\end{equation}
Where 
\[
F(a)= \int_{-\infty}^{\infty}f(\tau)e^{-a\tau}d\tau,
\]
which is a factor that depends on signal $f$. The convolution will be well defined only when $F(a)$ results in a finite value.
\begin{proof}
\begin{eqnarray*}
(f*g)(t) & = & \int_{-\infty}^{\infty}f(\tau)g(t-\tau)d\tau \\
         & = & \int_{-\infty}^{\infty}f(\tau)e^{a(t-\tau)}d\tau \\
         & = & \int_{-\infty}^{\infty}f(\tau)e^{at}e^{-a\tau}d\tau\\
         & = & \left[\int_{-\infty}^{\infty}f(\tau)e^{-a\tau}d\tau\right]e^{at}\\
         & = & F(a)g(t).
\end{eqnarray*}
\end{proof}
\item[(b) Discrete Convolution with exponential:] Let be $f$ a discrete signal and consider $g(k) = a^k$, with $a\neq 0\in\mathbb{C}$. 
Then 
\begin{equation}\label{convexpd}
(f*g)(k) = F(a)g(k)
\end{equation}
Where 
\[
F(a)= \sum_{n=-\infty}^{\infty}f(n)a^{-n},
\]
which is a factor that depends on signal $f$. The convolution will be well defined only when $F(a)$ is finite. 
\begin{proof}
\begin{eqnarray*}
(f*g)(k) & = & \sum_{n=-\infty}^{\infty}f(n)g(k-n)\\
         & = & \sum_{n=-\infty}^{\infty}f(n)a^{k-n} \\
         & = & \sum_{n=-\infty}^{\infty}f(n)a^ka^{-n}\\
         & = & \left[\sum_{n=-\infty}^{\infty}f(n)a^{-n}\right]a^k\\
         & = & F(a)g(k).
\end{eqnarray*}

\end{proof}

\end{description}

\end{proposition}

\begin{remark}\label{pconvexp}\em 
We also have an equivalent of Proposition~\ref{convexp} for periodic convolution:

\begin{description}
\item[(a) Analog Periodic Convolution with exponential:] Let be $f$ and $g$ analog periodic signals with period $T$ and consider 
$g$ the periodic signal obtained from the component $g_c(t) = e^{at}$ ($a\neq 0\in\mathbb{C}$) 
for $0\le t < T$ and zero otherwise). Then 
\begin{equation}\label{pconvexpa}
(f\circledast g)(t) = F(a)g(t)
\end{equation}
Where 
\[
F(a)= \int_{-T/2}^{T/2}f(\tau)e^{-a\tau}d\tau.
\]
\item[(a) Discrete Periodic Convolution with exponential:] Let be $f$ and $g$ discrete periodic signals with period $N$ and consider 
$g$  the periodic signal obtained from the component $g_c(k)=a^k$ ($a\neq 0\in\mathbb{C}$) for $0\le k \le N-1$ and 
zero otherwise. Then 
\begin{equation}\label{pconvexpd}
(f\circledast g)(k) = F(a)g(k)
\end{equation}
Where 
\[
F(a)= \sum_{n=0}^{N-1}f(n)a^{-n}.
\]

\end{description}

\end{remark}

Below we list some other properties of convolution that might be important for proving some properties of Fourier transform:

\begin{description}
\item[(i) Derivative of analog Convolution:] Let be $f$ and $g$ analog signals, with $f$ or $g$ differentiable 
(i.e. $\dot f$ or $\dot g$ exists):
\[
\dot f*g = f*\dot g = \dot{f*g}
\]

\item[(ii) Time shifting:] Let be $f$ and $g$ signals and we denote $[f]_a$ as the shifting of $f$ by ``$a$" units, that is:
$[f]_a(t) = f(t-a)$. Then:
\[
[f]_a*g = f*[g]_a = [f*g]_a
\]

\item[(iii) Time scaling:] let be $f$ signals and denote $f^a(t) = f(at)$ for $a\neq 0$, then:
\[
f^a*g = \frac{1}{|a|}(f*g^{1/a})^a \quad\text{or}\quad (f^a*g)(t) = \frac{1}{|a|}(f*g^{1/a})(at)
\]
Obs.: $g^{1/a}(t) = g(t/a)$

\end{description}

All properties also have their counterparts in discrete case. We note that, in fact, these properties show us how some  
operations can be ``transferred" from one signal to another under convolution.

\section{Fourier Series and Transforms as Convolution with Exponential}

\subsection{The Fourier Series}\label{sfourier}

It is well known that a analog periodic signal $f$ (with period $T$) can be written as an exponential Fourier series as shown 
below:\footnote{Of course there are some mathematical conditions that must be satisfied in order the Fourier series converges.
In particular, when $f$ is square integrable, over its period $T$, the series converges to $f(t)$ at almost every point 
$t\in\mathbb{R}$. Most signals used in engineering satisfies this condition of integrability.}
\begin{equation}\label{fseriesdef}
f(t) = \sum_{n=-\infty}^{\infty}C_ne^{jn\omega_0t},\quad \omega_0=2\pi/T
\end{equation}
and 
\begin{equation}\label{cndef}
C_n=\frac{1}{T}\int_{-T/2}^{T/2}f(\tau)e^{-jn\omega_0\tau}d\tau, \quad \omega_0=2\pi/T
\end{equation}
are the Fourier coefficients of the complex series. Also, if we consider $t\in\mathbb{R}$ representing time (e.g. seconds), we have that 
$\omega_0$ represents angular frequency (e.g radians/second), and Fourier coefficients $C_n$ may be seen as a (complex) discrete signal 
whose values are spaced by $\omega_0$ in frequency domain. 

To proceed with our analysis, we will first consider the complex exponential ``$e^{jn\omega_0t}$"as two different signals, as shown below:
\begin{description}
\item[(i) $n\in \mathbb{Z}$ is fixed:] $x_n(t) = e^{jn\omega_0t}$, $\omega_0=2\pi/T$, is a analog signal defined in time 
domain and $x_n$ is periodic since $x_n(t+T)=x_n(t)$, for all $t\in\mathbb{R}$.
\item[(ii) $t\in\mathbb{R}$ is fixed:] $x_t(n) = e^{jn\omega_0t}$, $\omega_0=2\pi/T$, is a discrete signal defined in 
frequency domain and whose values are spaced by $\omega_0$ (it is not necessarily periodic). 
\end{description}

We now present the main result, which corresponds to the Fourier series for a periodic signal:

\begin{proposition}\label{fsconv}\em 
Let it be a periodic analog signal $f$ with period $T$ and consider $x_n(t) = e^{jn\omega_0t}$ and $\bar x_t(n) = e^{-jn\omega_0t}$ with 
$\omega_0=2\pi/T$. Then we have the following pair of equations:
\begin{eqnarray}
(f\circledast x_n)(t) & = & F(n)x_n(t), \label{tform}\\
(F*\bar x_t)(n) & = & Tf(t)\bar x_t(n) \label{wform}
\end{eqnarray}

And 
\[
F(n) = \int_{-T/2}^{T/2}f(\tau)e^{-jn\omega_0\tau}d\tau = TC_n,
\]
where $C_n$ are the Fourier series coefficients of $f$ as defined in (\ref{cndef})
\begin{proof}
To prove Equation~(\ref{tform}) we use the fact that signals $f$ and $x_n$ are analog signals with same period $T$, 
and since $x_n(t) = e^{at}$ with $a=jn\omega_0$, the result is a consequence of 
convolution with exponential as shown in Remark~\ref{pconvexp}--Equation~(\ref{pconvexpa}): 
\[
(f\circledast x_n)(t) = F(a)x_n(t), \quad F(a) = \int_{-T/2}^{T/2}f(\tau)e^{-a\tau}d\tau, \quad a=jn\omega_0
\]
and we can represent $F(a)$ as $F(n)$, since $a=jn\omega_0$.


To prove Equation~(\ref{wform}), we note that $F(n)$ is a discrete signal (aperiodic 
in general), and its (discrete) convolution with (also discrete) exponential signal $\bar x_t(n)=a^n$, where $a=e^{-j\omega_0t}$, follows
directly from Proposition~\ref{convexp}--Equation~(\ref{convexpd}):

\[
(F*\bar x_t)(n)  = G(a)\bar x_t(n), \quad G(a) = \sum_{m=-\infty}^{\infty}F(m)a^{-m}, \quad a=e^{-j\omega_0t}
\]
So we have $G(a) = G(t)$, since $a=e^{-j\omega_0t}$, and 
\[
G(t) = \sum_{m=-\infty}^{\infty}F(m)e^{jm\omega_0t}, 
\]
We also have $F(m) = TC_m$, where $C_m$ are the Fourier coefficients of $f$,  then
\[
G(t) = \sum_{m=-\infty}^{\infty}TC_me^{jm\omega_0t} = T\sum_{m=-\infty}^{\infty}C_me^{jm\omega_0t} = Tf(t)
\]
and so we get 
\[
(F*\bar x_t)(n)  = G(a)\bar x_t(n) = G(t)\bar x_t(n) = Tf(t)\bar x_t(n).
\]
\end{proof}

%
\end{proposition}

\subsection{The Discrete Fourier Transform - DFT}\label{dtfourier}
Discrete Fourier transform or DFT is a version of Fourier series when signal $f$ is discrete with period $N>0\in\mathbb{Z}$:

\begin{proposition}\label{dftconv}\em 
Let it be a periodic discrete signal $f$ with period $N$ and consider $x_n(k) = e^{jn(2\pi/N)k}$ and $\bar x_k(n) = e^{-jk(2\pi/N)n}$ 
also both periodic with period $N$. Then we have the following pair of equations:
\begin{eqnarray}
(f\circledast x_n)(k) & = & F(n)x_n(k) \label{tform3}\\
(F\circledast\bar x_k)(n) & = & Nf(k)\bar x_k(n) \label{wform3}
\end{eqnarray}
And 
\[
F(n) = \sum_{m=0}^{N-1}f(m)e^{-jm(2\pi/N)n}
\]
which is periodic with period $N$, since $F(n+N) = F(n)$ for all $n$.
$F$ is denominated Discrete Fourier Transform (or DFT) of $f$.

\begin{proof}
To prove Equation~(\ref{tform3}) we use the fact that signals $f$ and $x_n$ are analog signals with same period $N$, 
and since $x_n(k) = a^k$ with $a=e^{jn2\pi/N}$, the result is a consequence of 
convolution with exponential as shown in Remark~\ref{pconvexp}--Equation~(\ref{pconvexpd})
\[
(f\circledast x_n)(k) = F(a)x_n(k), \quad F(a) = \sum_{m=0}^{N-1}f(m)a^{-m}, \quad a=e^{jn2\pi/N}
\]
and we can represent $F(a)$ as $F(n)$.

Before proceeding to prove (\ref{wform3}) we use (\ref{tform3}) to prove the following ``ortogonality" 
condition between periodic exponential discrete signals 
$x_m(k) = e^{jm(2\pi/N)k}$ and $x_n(k) = e^{jn(2\pi/N)k}$:
\begin{corollary}\label{dxmxn}\em 
Let it be the periodic signals $x_m(k) = e^{jm(2\pi/N)k}$ and $x_n(k) = e^{jn(2\pi/N)k}$, then:
\[
(x_m\circledast x_n)(k) = N\delta(m-n)x_n(k), \quad\text{with } 
\delta(m-n)=\begin{cases}1, & \text{if } m=n \\ 0, & \text{otherwise}\end{cases}
\]
And so, we have $(x_n\circledast x_n)=Nx_n$ and $(x_m\circledast x_n)=0$ for $m\neq n$.
\begin{proof}
Signals $x_m$ and $x_n$ have same period $N$ and then considering $f=x_m$ in Equation~(\ref{tform3}) we easily get $F(n) = N\delta(m-n)$ 
by solving the summand.
\end{proof}
\end{corollary}

We now proceed to prove Equation~(\ref{wform3}). We have that $F$ and $\bar x_k$ are both periodic with same period $N$, and 
since we can write $\bar x_k(n) = a^n$, with $a=^{-jk2\pi/N}$ we again use the result of convolution with exponential as 
shown in Remark~\ref{pconvexp}--Equation~(\ref{pconvexpd}):
\[
(F\circledast\bar x_k)(n) = G(a)\bar x_k(n), \quad G(a) = \sum_{m=0}^{N-1}F(m)a^{-m}, \quad a=^{-jk2\pi/N}
\]
and we can represent $G(a)$ as $G(k)$. 
We will show that, in fact, $G(k)=Nf(k)$, and for that we use the ``ortogonality" result of Corollary~\ref{dxmxn}:
\begin{eqnarray*}
G(k) & = & \sum_{m=0}^{N-1}F(m) x_m(k)\\
(G\circledast x_n)(k) & = & \sum_{m=0}^{N-1}F(m) (x_m\circledast x_n)(k)\\
                      & = & \sum_{m=0}^{N-1}F(m)(N\delta(m-n))x_n(k)\\
                      & = & NF(n)x_n(k) = N(f\circledast x_n)(k), \quad\text{by } (\ref{tform3})
\end{eqnarray*}
And so we have
\[
(G\circledast x_n)(k) = (Nf\circledast x_n)(k) \implies G(k) = Nf(k),
\]
which can be easily shown by solving a simple non-singular linear system with $N$ equations and $N$ unknowns.

\end{proof}

\end{proposition}

\subsection{The Fourier Transform}\label{tfourier}

We will present the Fourier transform as a limit case of the Fourier series, as shown in Proposition~\ref{fsconv}, when period $T$ of signal 
$f$ tends to infinity.

\begin{proposition}\label{ftconv}\em 
Let be $f$ an absolutely integrable analog signal and consider the analog signals $x_{\omega}(t) = e^{j\omega t}$ and 
$\bar x_t(\omega)=e^{-j\omega t}$, then we have the following pair of equations:
\begin{eqnarray}
(f*x_{\omega})(t) & = & F(\omega)x_{\omega}(t) \label{tform1}\\
(F*\bar x_t)(\omega) & = & 2\pi f(t)\bar x_t(\omega) \label{wform1}
\end{eqnarray}
And
\[
F(\omega) = \int_{-\infty}^{\infty}f(\tau)e^{-j\omega\tau}d\tau 
\]
is the Fourier Transform of $f$.

\begin{proof}
We consider initially $f$ as being a periodic signal with period $T=2\pi/\omega_0$ and so, by Proposition~\ref{fsconv}, 
we have the following pair
\begin{eqnarray*}
(f\circledast x_n)(t) & = & F(n)x_n(t) \\
(F*\bar{x}_t)(n) & = &Tf(t)\bar{x}_t(n)=\frac{2\pi}{\omega_0}f(t)\bar{x}_t(n)
\end{eqnarray*}
Equivalently
\begin{eqnarray}
(f\circledast x_n)(t) & = & F(n)x_n(t) \label{tform2} \\
\omega_0(F*\bar{x}_t)(n) & = &2\pi f(t)\bar{x}_t(n) \label{wform2}
\end{eqnarray}
Now we make $T\to\infty$ and so $\omega_0 \to 0$ which it is an infinitesimal ``$d\omega$". Similarly 
we have done before in Definition~\ref{aconv}, when $\omega_0=d\omega$ we have $n\omega_0\to \omega$, $F(n)\to F(\omega)$,
$\bar x_t(n)\to \bar x_t(\omega)$, since $\omega_0$ is the spacing of the values of $F(n)$ (and also of $\bar x_t(n)$) in frequency domain. 
Then the discrete convolution in left-hand side of Equation~(\ref{wform2}) 
turns into an analog convolution between $F(\omega)$ and $\bar x_t(\omega)$. On the other hand, 
the circular analog convolution in left hand side of Equation~(\ref{tform2}) turns into a (regular) analog convolution when $T\to\infty$. So we get the pair of Equations (\ref{tform1}) and
(\ref{wform1}). Finally, we note that Equation~(\ref{tform1}) is essentially Equation~(\ref{convexpa}) in Proposition~\ref{convexp} 
(with $a=j\omega$) and so
\[
F(\omega) = \int_{-\infty}^{\infty}f(\tau)e^{-j\omega\tau}d\tau,
\]
which is the Fourier Transform of $f$.
\end{proof}

\end{proposition}

\section{Applications}\label{aplicat}

\subsection{Fourier Series}\label{fsapp}

The formulation of Fourier series presented in Proposition~\ref{fsconv}, in our view, simplify proofs for some Fourier series
properties. We list some of them below:

\begin{description}
\item[(a) Convolution in time:] Let be $f$ and $g$ periodic (with same period). Which is the spectrum of their circular convolution?
\[
(f\circledast x_n)(t) = F(n)x_n(t), \qquad (g\circledast x_n)(t) = G(n)x_n(t)
\]
Then
\begin{eqnarray*}
[(f\circledast g)\circledast x_n](t) & = & [f\circledast(g\circledast x_n)](t)\\
                                     & = & [f\circledast(G(n)x_n)](t) \\
                                     & = & G(n)(f\circledast x_n)(t) \\
                                     & = & [G(n)F(n)]x_n(t)
\end{eqnarray*}

\item[(b) Convolution in frequency:] Which periodic signal is obtained by the (discrete) convolution between the spectra of $f$ and $g$, 
which are periodic with same period?
\[
(F*\bar x_t)(n) = Tf(t)\bar x_t(n), \qquad (G*\bar x_t)(n) = Tg(t)\bar x_t(n)
\] 
Then
\begin{eqnarray*}
[(F*G)*\bar x_t](n) & = & [F*(G*\bar x_t)](n)\\
                    & = & [F*(Tg(t)\bar x_t)](n) \\
                    & = & Tg(t)(F*\bar x_t)(n) \\
                    & = & Tg(t)Tf(t)\bar x_t(n) \\
                    & = & T[Tg(t)f(t)]\bar x_t(n)
\end{eqnarray*}

\item[(c) Convolution in time with an aperiodic signal:] Let be $h$ an aperiodic (and absolutely integrable) signal and $u$ a periodic signal. 
Which is the spectrum of the periodic signal ``$h*u$"?
\[
(u\circledast x_n)(t) = U(n)x_n(t), \quad (h*x_n)(t) = H(n)x_n(t)
\]
We note that ``$H(n)$" exists since ``$h$" is absolutely integrable. Then

\begin{eqnarray*}
[(h*u)\circledast x_n](t) & = & [h*(u\circledast x_n)](t) \\
                          & = & [h*(U(n)x_n)](t)\\
                          & = & U(n)(h*x_n)(t)\\
                          & = & [U(n)H(n)]x_n(t)
\end{eqnarray*}
Obs.: We can see ``$U(n)H(n)$" as the spectrum of the output signal of a 
{\em stable Linear and Time-Invariant system} with impulse response ``$h$", when the input is a periodic signal ``$u$".
\end{description}

We believe other properties can be easily deduced from the formulation proposed in Proposition~\ref{fsconv} for the Fourier series.

\subsection{Fourier Transforms}\label{tfourierapp}

We will derive some properties of Fourier transforms using the formulation presented in Proposition~\ref{ftconv}.

\begin{description}

\item[(a) Convolution in time:] Let be $f$ and $g$ with Fourier transform $F$ and $G$, respectively. 
Which is the Fourier transform of $f*g$?
\[
(f*x_{\omega})(t) = F(\omega)x_{\omega}(t), \qquad (g*x_{\omega})(t) = G(\omega)x_{\omega}(t)
\]
Then
\begin{eqnarray*}
[(f*g)*x_{\omega}](t) & = & [(f*(g*x_{\omega})](t) \\
                    & = & [f*(G(\omega)x_{\omega}](t) \\
                    & = & G(\omega)(f*x_{\omega})(t) \\
                    & = & [G(\omega)F(\omega)]x_{\omega}(t)
\end{eqnarray*}

\item[(a) Convolution in Frequency:] Let be $f$ and $g$ with Fourier transform $F$ and $G$, respectively. 
Which is the inverse Fourier transform of $F*G$?
\[
(F*\bar x_t)(\omega) = 2\pi f(t)\bar x_t(\omega), \qquad (G*\bar x_t)(\omega)=2\pi g(t)\bar x_t(\omega)
\]
Repeating the reasoning used before in item (a), we easily obtain
\[
[(F*G)*\bar x_t](\omega) = 2\pi[2\pi f(t)g(t)]\bar x_t(\omega)
\]

\item[(c) Derivative in time:] Given the Fourier transform of $f$ (differentiable) obtain (when exists) the Fourier transform of $\dot f$.
\[
(f*x_{\omega})(t) = F(\omega)x_{\omega}(t), \quad x_{\omega}(t) = e^{j\omega t}
\]
Then
\[
(\dot f*x_{\omega})(t) = (f*\dot x_{\omega})(t) = [f*(j\omega x_{\omega})](t) = j\omega(f*x_{\omega})(t) = 
j\omega F(\omega)x_{\omega}(t).
\]

\item[(d) Shifting in time:] Let $f$ with Fourier transform $F$. which is the Fourier transform for $[f]_{t_0}(t) = f(t-t_0)$?
\[
(f*x_{\omega})(t) = F(\omega)x_{\omega}(t), \quad x_{\omega}(t) = e^{j\omega t}
\]
Then
\[
([f]_{t_0}*x_{\omega})(t) = (f*[x_{\omega}]_{t_0})(t) = [f*(e^{-j\omega t_0}x_{\omega})](t) = e^{-j\omega t_0}(f*x_{\omega})(t) = 
e^{-j\omega t_0}F(\omega)x_{\omega}(t)
\]

\item[(e) Duality:] Let be $f(t)$ with Fourier transform $F(\omega$). Which is the Fourier transform of $F(t)$?
\[
(F*x_\omega)(t) = G(\omega)x_{\omega}(t), \quad\text{who is } G(\omega)\;?
\]
We have
\begin{eqnarray*}
(F*\bar x_t)(\omega) & = & 2\pi f(t)\bar x_t(\omega), \quad t \leftrightarrows \omega\\
(F*\bar x_{\omega})(t) & = & 2\pi f(\omega)\bar x_{\omega}(t), \quad \omega \to -\omega \\
(F*\bar x_{-\omega})(t) & = & 2\pi f(-\omega)\bar x_{-\omega}(t), \quad \bar x_{-\omega}(t) = x_{\omega}(t) \\
(F*x_{\omega})(t) & = & \underbrace{[2\pi f(-\omega)]}_{G(\omega)}x_{\omega}(t)
\end{eqnarray*}

\item[(f) Time scaling:] Let be $f$ with Fourier transform $F$. Which the Fourier transform of $f^a$, where $f^a(t) = f(at)$?
\[
(f*x_{\omega})(t) = F(\omega)x_{\omega}(t), \quad x_{\omega}(t) = e^{j\omega t}
\]
Then
\begin{eqnarray*} 
(f^a*x_{\omega})(t) & = & \frac{1}{|a|}(f*x_{\omega}^{1/a})(at), \quad x_{\omega}^{1/a} = x_{\omega/a}\\
                    & = & \frac{1}{|a|}(f*x_{\omega/a})(at)\\
                    & = & \frac{1}{|a|}F(\omega/a)x_{\omega/a}(at), \quad x_{\omega/a}(at)= x_{\omega}(t) \\
                    & = & \frac{1}{|a|}F(\omega/a)x_{\omega}(t)
\end{eqnarray*}

\item[(g) Discretization of Fourier Transform:] Let be $f(t)$ with Fourier Transform $F(\omega)$. How do we
interpret $F(n\omega_0)$, when $\omega_0$ is a interval in frequency domain?

We have
\begin{eqnarray}
(f*x_{\omega})(t) & = & F(\omega)x_{\omega}(t), \quad x_{\omega}(t) = e^{j\omega t} \label{tform41}\\
(F*\bar x_t)(\omega) & = & 2\pi f(t)x_t(\omega), \quad \bar x_t(\omega) = e^{-j\omega t} \label{wform41}
\end{eqnarray}

Taking $\omega = n\omega_0$ we have 
\begin{eqnarray*}
x_{\omega}(t) = x_{n}(t) & = & e^{jn\omega_0t} \quad\text{(become periodic with period } T=2\pi/\omega_0)\\
\bar x_t(\omega) = \bar x_t(n) & = & e^{-jn\omega_0t} \quad\text{(become discrete})
\end{eqnarray*}
In equation~(\ref{tform41}), changing $\omega$ by $n\omega_0$ we get:
\begin{equation}\label{tformperiod}
(f * x_n)(t) = F(n)x_n(t), \quad x_n(t) = e^{jn\omega_0t}
\end{equation}
and, since $x_n(t)$ is periodic with period $T=2\pi/\omega_0$, we have that the convolution in left-hand side of 
(\ref{tformperiod}) is periodic. 
We now define (from $f$) a signal
$f_p$ periodic also with period $T$ so that $f$ corresponds to $f_p$ over one period.\footnote{ A condition may be imposed on $f$ 
such that this works properly: in fact $f$ must ``fit" in interval $T$ and so must be non-null only over a finite interval of 
time to avoid the so-called aliasing phenomenon.} Then, as discussed in 
Definition~\ref{pconvdef}, the (regular) convolution ($f*x_n$) in left-hand side of (\ref{tformperiod}) can be 
re-written as a circular convolution $(f_p\circledast x_n)$, and then:
\[
(f_p\circledast x_n)(t) = F(n)x_n(t)
\]
By Proposition~\ref{fsconv}, we have that $F(n)=TC_n$, where $C_n$ are the Fourier coefficients of the series of periodic signal
$f_p$. 

\item[(h) Spectrum of a sampled signal:] Let be $f(t)$ with Fourier Transform $F(\omega)$. Analyze the spectrum of $f^*(k) = f(kT_s)$. 
Let be the Fourier transform pair of $f$:
\begin{eqnarray}
(f*x_{\omega})(t) & = & F(\omega)x_{\omega}(t), \quad x_{\omega}(t) = e^{j\omega t} \label{tform4}\\
(F*\bar x_t)(\omega) & = & 2\pi f(t)x_t(\omega), \quad \bar x_t(\omega) = e^{-j\omega t} \label{wform4}
\end{eqnarray}
and take $t = kT_s$, so that we have:
\begin{eqnarray*}
x_{\omega}(t) = x_{\omega}(k) & = & e^{jkT_s\omega} \quad\text{(become discrete)}\\
\bar x_t(\omega) = \bar x_k(\omega) & = & e^{-jkT_s\omega} \quad\text{(become periodic with period } \omega_s=2\pi/T_s)
\end{eqnarray*}

Additionally, by making $t=kT_s$ in Equation~(\ref{wform4}) we get:
\begin{equation}\label{wformperiod}
(F *\bar x_k)(\omega) = 2\pi f^*(k)\bar x_k(\omega) 
\end{equation}
And we note that convolution in left-hand side of Equation~(\ref{wformperiod}) is now periodic, since
$\bar x_k(\omega)$ is periodic with period $\omega_s=2\pi/T_s$. We define a signal $F^*(\omega)$ 
(from $F$, the Fourier transform of $f$) that is periodic with period $\omega_s$ so that $F$ corresponds to $F^*$ 
over one period.\footnote{ A similar condition, as 
discussed in (g), have to be imposed now on $F$.} Then, as discussed in 
Definition~\ref{pconvdef}, the (regular) convolution ($F*\bar{x}_k$) in left-hand side of (\ref{wformperiod}) can be 
re-written as a circular convolution:
\[
(F^* \circledast\bar x_k)(\omega) = 2\pi f^*(k)\bar x_k(\omega) 
\quad\text{ or }\quad (F^*/(2\pi) \circledast\bar x_k)(\omega) = f^*(k)\bar x_k(\omega) 
\]
In order to obtain a standard format of Fourier series equations as shown in (Proposition~\ref{fsconv}), 
we make $k\to -k$ so that $\bar x_k = x_k$ and 
and define $G(\omega) = F^*(\omega)/(2\pi)$ and also $g(k) = f(-k)$. Then equation 
above can be re-written as: 

\begin{equation}\label{wformcirc}
(G \circledast x_k)(\omega) = g(k)x_k(\omega).
\end{equation}

By comparing Equation~(\ref{wformcirc}) above with Equation~(\ref{tform}), we note that
$g(k)$ can be seen as the ``spectrum" of periodic signal $G$ obtained by turning $F$ (Fourier transform of $f$)
periodic with period $\omega_s=2\pi/T_s$. By result of Proposition~\ref{fsconv} we have the following pair of 
equations in Fourier series format:
\begin{eqnarray}
(G \circledast x_k)(\omega) & = & g(k)x_k(\omega) \label{wtform}\\
(g*\bar x_{\omega})(k) & = & \omega_sG(\omega)\bar{x}_{\omega}(k) \label{twform}
\end{eqnarray}

With $g(k) = f^*(-k)$ and $G(\omega) = F^*(\omega)/(2\pi)$, where $F^*(\omega)$ is a periodic signal, with period 
$\omega_s=2\pi/T_s$ ($T_s$ is the sampling time interval), such that one period of it
corresponds to $F$, the Fourier Transform of $f$.

\end{description}

\subsection{DFT}\label{dftapp}

\begin{description}
\item[(a) DFT versus Fourier series:] 
Lets suppose we have $N$ samples of a (periodic) signal $f$, which are supposed to be obtained from one period $T$ of $f$, and 
additionally they are spaced in time by a sampling interval $T_s$ (so that $T/T_s=N$). 
We represent these samples by discrete signal $f_d$ and its DFT by $F_d(n)$. Then by Proposition~\ref{dftconv} we have
\begin{equation}\label{dftfdt}
(f_{dt}\circledast x_n)(k) = F_d(n)x_n(k) 
\end{equation}
Where $f_{dt}$ represents a periodic signal obtained by repeating $f_d$. 

With samples $f_d$ we can obtain the Fourier transform $F_a(\omega)$ of analog signal $f_a$ (one period of $f$), 
whose samples results in $f_d$, by
using Equation~(\ref{twform}):
\begin{equation} \label{twform1}
(g*\bar x_{\omega})(k)  = \omega_sG(\omega)\bar{x}_{\omega}(k)
\end{equation}
where $g(k) = f_d(-k)$ and $G(\omega) = F_a^*(\omega)/(2\pi)$. We have that $F_a^*$ is a periodic repetition of $F_a$ with
period $\omega_s=2\pi/T_s$. To make $f_d$ appear in (\ref{twform1}) we make $k\to -k$ and so we have:
\begin{equation}\label{twform2}
(f_d*x_{\omega})(k) = \frac{\omega_s}{2\pi}F_a^*(\omega)x_{\omega}(k) =
\frac{F_a^*(\omega)}{T_s}x_{\omega}(k)
\end{equation}
Lets consider $\omega=n\omega_0$, with $\omega_0=2\pi/T$ and represent $x_{\omega}$ as $x_n$ and $F_a^*(\omega)$ as $F_a^*(n)$ in 
Equation~(\ref{twform2}) above to obtain:
\begin{equation}\label{fttemp}
(f_d*x_{n})(k) = \frac{F_a^*(n)}{T_s}x_{n}(k)
\end{equation}
Since $f_d$ is one period of $f_{dt}$, the convolution in left-hand side of (\ref{fttemp}) is periodic, and we have
\begin{equation}\label{fttemp1}
(f_{dt}\circledast x_{n})(k) = \frac{F_a(n)}{T_s}x_{n}(k)
\end{equation}
and $F_a(n)$ is one period of $F_a^*(n)$. Since $F_a(n)$ is one period of $F_a^*(n)$, then $F_a(n)$ is the discretization 
of $F_a(\omega)$, i.e., the Fourier transform of one period of 
$f$ (namely $f_a$), and so $F_a(n) = TC_n$, where $C_n$ are the Fourier coefficients of periodic signal $f$ as shown in item (g) 
of Section~\ref{tfourierapp}. We then re-write (\ref{fttemp1}) as 
\begin{equation}\label{fttemp2}
(f_{dt}\circledast x_{n})(k) = \frac{TC_n}{T_s}x_{n}(k) = (NC_n)x_n(k)
\end{equation}
Comparing Equation~(\ref{fttemp2}) with Equation~(\ref{dftfdt}), we get $F_d(n)=NC_n$, where $F_d$ if the discrete Fourier transform
of $f_{dt}$ and $C_n$ are the Fourier coefficients of analog signal $f$ whose discretization ($N$ samples by period) results in $f_{dt}$.

\end{description}

\section{Conclusions}

We have shown in this note that the Fourier Series and Transform can be formulated as a set of two equations involving 
a convolution with an exponential signal, where in one of the equations the frequency is fixed and in another the time is fixed. 
We used the idea to show how to prove some properties of Fourier series and the Fourier transform, and given its simplicity,
we think it could be useful as an alternative approach for the study of Fourier Series and transforms. We also mention that other 
transforms, like Laplace and $Z$, also can be formulated in this way and may could be interesting to be analyzed.

\bibliographystyle{plain}

\end{document}